\documentclass{article}
\usepackage{epubtk}   
\usepackage{booktabs}
\usepackage{graphicx}
\usepackage{graphics}
\usepackage{lscape}
\usepackage{indentfirst}
\usepackage{latexsym}
\usepackage{authblk}
\usepackage{blindtext}
\usepackage{amsmath}
\usepackage{amssymb}
\usepackage{amsthm}
\usepackage{algorithm}
\usepackage{algpseudocode} 
\usepackage{color}
\usepackage{amssymb}
\usepackage{amsfonts}
\usepackage{amsmath}
\usepackage{bm}
\usepackage{mathrsfs}

\newcommand{\bea}{\begin{eqnarray}}
\newcommand{\eea}{\end{eqnarray}}

\newcommand{\cI}{{\cal I}}
\newcommand{\be}{\begin{equation}}
\newcommand{\ee}{\end{equation}}
\newcommand{\blank}{\bigskip}
\newcommand{\lam}{\lambda}
\newcommand{\La}{\Lambda}
\newcommand{\ka}{\kappa}
\newcommand{\ft}{\footnote}

\newtheorem*{theorem}{Theorem}

\title{ On the stability and accuracy of the Empirical Interpolation Method and Gravitational Wave Surrogates}
\author{Manuel Tiglio  \\ Aar\'on Villanueva}
\affil{Facultad de Matem\'atica, Astronom\'ia, F\'isica y Computaci\'on, Universidad Nacional de C\'ordoba (5000), C\'ordoba, Argentina}

\begin{document}

\maketitle

\begin{abstract}
The combination of the Reduced Basis and the Empirical Interpolation Method (EIM) approaches have produced outstanding results in many disciplines. In particular, in gravitational wave (GW) science these results range from building non-intrusive surrogate models for GWs to fast parameter estimation adding the use of Reduced Order Quadratures. These surrogates have the salient feature of being essentially indistinguishable from or very close to supercomputer simulations of the Einstein equations, but can be evaluated in the order of a milisecond per multipole mode on a standard laptop. In this article we clarify a common misperception of the EIM as originally introduced and used in practice in GW science. Namely, we prove that the EIM at each iteration chooses the interpolation nodes so as to make the related Vandermonde-type matrix as invertible as possible; not necessarily optimizing its conditioning or accuracy of the interpolant as is sometimes thought. In fact, we introduce two new variations of the EIM, nested as well,  which do optimize with respect to conditioning and the Lebesgue constant, respectively, and compare them through numerical experiments with the original EIM using GWs. Our analyses and numerical results suggest a subtle relationship between solving for the original EIM, conditioning, and the Lebesgue constant, in consonance with active research in rigorous approximation theory and related fields. 
\end{abstract}
%------------------------------------------------------------------------------------------------
\section{Introduction} 
%------------------------------------------------------------------------------------------------
The general problem of determining existence, uniqueness, stability, and accuracy of an optimal interpolant for an arbitrary basis, including the choice of interpolation nodes, is still largely an open one; for recent results in the least squares sense see \cite{Cohen2013OnTS,CohenCorr,Cohen2016} . The one-dimensional polynomial case is well understood. For example, it is well known that the Chebyshev nodes constitute an optimal solution as they satisfy a min-max property (see \cite{cheb59} or \cite{Deuflhard2003NumericalAI}, pag. 195). Yet, Gaussian nodes are not hierarchical (nested), which is in many cases a desired condition.

The multi-dimensional case is more complicated and less developed, as the evaluation of tensor products of one-dimensional polynomials typically suffer from the curse of dimensionality \cite{NarayanX12} and are also restricted to simple geometries. An active area of research is that one of sparse grids \cite{Garcke2013}; however, most of the approaches are not nested and in many cases not applicable to complex geometries.

Over the last decade a general quasi-optimal approach for interpolation of parameterized problems, which works extremely well in practice, has gained wide popularity: the Empirical Interpolation Method (EIM)\cite{Barrault2004667, sorensen2010}, closely related to the Magic Points approach \cite{Maday_2009}. Given a reduced basis distilled from a parametrized fiducial model, the EIM is, in its original formulation, an efficient algorithm for finding a good set of nodal points and the corresponding interpolant. As of this writing, in GW science the combination of the EIM with a Proper Orthogonal Decomposition (POD) or a greedy approach to construct reduced bases constitute a powerful way of building accurate surrogate models and reduced order quadratures for binary systems without any physical approximation; for a recent review see Ref. \cite{Tiglio:2021ysj}. 

In this article we discuss some fundamental aspects of the EIM as currently used when building GW surrogates; namely, its existence, stability, and accuracy. These features have been observed and qualitatively discussed in the literature of GW surrogates; here we analyze them in more detail. Most of our analyses are actually independent of our application domain of interest and might be useful in other contexts. The approach here adopted is data-driven, that is, no differential equations are invoked. The intended audience of this work is that one of practitioners in GW science, not experts in approximation theory. The style here used aims to bring closer to the community some results from the backend mathematical theory that might not be so obviously seizable.

The structure of this article is as follows. In Section~\ref{sec:interp} we briefly review the general interpolation problem, and motivate the three main points that we analyze in this paper within the context of the EIM; namely: existence, conditioning and accuracy. In Section~\ref{sec:eim} we discuss the particular case of the EIM approach, as originally introduced. In Section \ref{sec:exist} we prove our main result, showing that the EIM at each iteration chooses the next interpolation node as the one which makes the determinant of the Vandermonde related matrix as large as possible. To our knowledge, this is a new result in the literature and is in contradiction, at least within our domain of application, to the common belief that the EIM optimizes for accuracy (however, as we will point out, there seems to be a subtle relationship between these two optimization criteria). Next, motivated by this result, in Sections \ref{sec:cond} and \ref{sec:acc} we discuss the conditioning and accuracy of the EIM, and in Section \ref{sec:opt} we present two new variations of the EIM which instead of optimizing with respect to the invertibility of the associated Vandermonde matrix optimize with respect to its conditioning and Lebesgue constant of the interpolant, respectively. Again, to our knowledge these new variations and their analyses here presented are new in the literature.  In Section~\ref{sec:numerics} we present numerical results from gravitational waves emitted by black hole collisions, initially without spin and in quasi-circular orbit, and compare the original EIM with respect to our variations.  We close in Section~\ref{sec:comments} with a summary of our results and some comments. 

%------------------------------------------------------------------------------------------------
\section{Interpolation} \label{sec:interp}
%------------------------------------------------------------------------------------------------

Since our motivation is that one of gravitational waves, we focus on time series --the frequency domain case is identical. 
The setting of the problem is as follows. GWs are commonly represented as complex functions of time indexed by a continuous --in general, multi-dimensional-- parameter $\lam$:

$$
h_{\lam}(t) = h(\lam, t) \,.
$$
The parameter $\lam$ might comprise the mass and the spin of a binary component, for example. The total space of GWs for a given physical setting compose a parametrized model. The goal is to represent it by a set of sparse basis functions $\{e_i(t)\}_{i=1}^n$ that leads to a surrogate, predictive model of high accuracy which is fast to evaluate. Ideally, one wants the surrogate to be indistinguishable from numerical relativity simulations of the full Einstein equations. The construction of such surrogates not only requieres a (ideally, sparse) basis, but, as of date, also interpolation in parameter space and time. 

The general interpolation problem consists of finding an interpolant of the form
\be
\cI_n[h_{\lam}](t)=\sum_{i=1}^n C_i(\lam) e_i(t) \, ,  \label{eq:interp1}
\ee
subject to the conditions  
\be
\cI_n[h_{\lam}](T_j)=h_{\lam}(T_j)\,,\,j=1,\ldots,n \label{eq:interp2}
\ee
for given nodes $\{ T_j \}$. 

The following are in general minimum requirements for any interpolation approach: 
\begin{enumerate}
\item Existence of the interpolant.
\item The solutions to the problems of defining/finding the location of the interpolation nodes, as well as constructing and evaluating the interpolant, have to be stable/well conditioned.
\item The accuracy, when compared to a projection-based approach, has to be competitive.
\end{enumerate}

Some extra conditions to ask for, depending on the application, might include:

\begin{enumerate}
\setcounter{enumi}{3}
\item The interpolation nodes are hierarchical (nested). 
\item The interpolant is fast to evaluate. 
\end{enumerate}
In this article we discuss items (1-3) in detail, take (4) for granted, and we briefly discuss (5), all within the context of the EIM.
 
Within the Reduced Basis Method (RBM) framework~\cite{maday2006,Fares20115532,hesthaven2015certified,quarteroni2015reduced}, the functions $\{ e_i \}$ in Eq.~(\ref{eq:interp1}) constitute a set of $n$ basis elements, called the {\it reduced basis}. It is found by an earlier implementation of, for example, POD or greedy-type algorithms. The span of the reduced basis gives a linear function that approximates the functions of interest up to a an arbitrary tolerance. We use the basis to define the $n\times n$ Vandermonde-like matrix (the ``V-matrix'') as
$$
 ({\bf V}_n)_{ij}:=e_j(T_i)\,,
 $$
where the label $n$ is to keep track of the basis dimension. Then, the interpolation conditions~(\ref{eq:interp2}) can be written as a linear system of the form 
 \be
 \sum_{i=1}^n ({\bf V}_n)_{ji}C_i (\lam) = h_{\lam}(T_j)\,, \quad j=1,\ldots,n  \label{eq:interp-problem}
 \ee
for the $n$ coefficients $\{ C_i \}$.

Notice that solving this interpolation problem does {\em not} involve any training set --which is defined as the discrete set of functions used to build the reduced basis-- but only knowledge of the basis. This is a highly non-trivial complexity reduction, since if the reduced basis is sparse the interpolation problem is a small one to solve for. This is in contrast to the construction of the basis itself, which can be computationally very expensive (see, for example, \cite{random-training}). In addition, it is also non-trivial that by using only the basis for the construction of the interpolant, the latter has a ``good'' representation error for {\em any} function in the space of interest. The intuition is that if the reduced basis is an accurate representation of this space, by transitivity so is the interpolant --assuming it was constructed with the choice of ``good enough'' interpolating nodes. One of the purposes of this article is to analyze this point in detail; namely, the accuracy of the EIM compared to variations of it that we here introduce.

If the V-matrix is invertible the solution to the interpolation problem (\ref{eq:interp1},\ref{eq:interp2}) is well defined and can be written in a compact manner:
\be
\cI_n[h_{\lam}](t) = \sum_{i=1}^n B_i(t) h_{\lam} (T_i) \, , \label{eq:interp-sol1}
\ee
where
\be
B_i (t) := \sum_{j=1}^n ({\bf V}_n^{-1})_{ji}e_j (t) \,.  \label{eq:interp-sol2}
\ee
The functions $B_i(t)$ are generalizations of the Lagrange functions in polynomial interpolation and satisfy $B_i(T_j)=\delta_{ij}$ for $i,j=1,\dots,n$.

For the purposes of this article, then, the items of required or desired properties of the interpolant above listed can be rephrased as:
\begin{enumerate}
\item Ensuring that the V-matrix is invertible.
\item Ensuring that it is not only invertible but also well conditioned, so that the inversion does not amplify numerical errors.
\item Analyzing the accuracy of the resulting interpolant.
\end{enumerate}
Items (1-3) are discussed in Sections \ref{sec:exist} , \ref{sec:cond} and \ref{sec:acc}, respectively. We briefly discuss computational cost in Section~\ref{sec:comments}. 

%------------------------------------------------------------------------------------------------
\section{The Empirical Interpolation Method (EIM)} \label{sec:eim}
%------------------------------------------------------------------------------------------------
An efficient way of selecting the nodes $\{T_j\}_{j=1}^n$ is by means of the EIM. As reviewed below, the algorithm receives as input the reduced basis and, by a hierarchical approach to the problem (\ref{eq:interp-problem}), it finds the empirical nodes and the associated interpolant. The first formulations of the algorithm was introduced and studied in \cite{Barrault2004667,Maday_2009} and its discrete implementation as is used for GWs, in \cite{sorensen2010}. See for example, \cite{PhysRevX.4.031006}.

{\scriptsize
\begin{algorithm}[H]
\caption{The Empirical Interpolation Method}
\label{alg:EIM}
\begin{algorithmic}[1]
\State {\bf Input:} $\{ e_i \}_{i=1}^n$
\vskip 10pt
\State $T_1 = \text{argmax}_t | e_1|$ 
\For{$j = 2 \to n$} 
\State Build ${\cal I}_{j-1} [e_j](t)$
\State $r_j(t) = e_j(t)-{\cal I}_{j-1} [e_j](t)$  (where $r$ stands for the {\it residual})
\State $T_j = \text{argmax}_t |r_j|$
 \EndFor
\vskip 10pt
\State {\bf Output:} EIM nodes $\{ T_i \}_{i=1}^n$ and interpolant $\cI_n$
\end{algorithmic}
\end{algorithm}
}

%------------------------------------------------------------------------------------------------
\subsection{Existence} \label{sec:exist} 
%------------------------------------------------------------------------------------------
We start by analyzing existence of the EIM at each iteration, which leads to our main result in the form of a theorem. We begin with an observation about the 5th step of Algorithm \ref{alg:EIM}. The first iteration implies the construction of the interpolant (Step 4)
$$
{\cal I}_{1} [e_2](t)=\frac{e_2(T_1)}{e_1(T_1)}e_1(t) \, , 
$$
and the associated residual (Step 5)
\be\label{eq:r_2}
r_2(t)=e_2(t)-\frac{e_2(T_1)}{e_1(T_1)}e_1(t)\,.
\ee
If we rewrite (\ref{eq:r_2}) as 
\be\label{eq:e2}
r_2(t)=\frac{e_2(t)e_1(T_1)-e_2(T_1)e_1(t)}{e_1(T_1)}\,
\ee
one recognizes in the numerator the determinant of
$$
\label{eq:V_2}
  {\bf V}_2(T_1, t) := \left(  \begin{array}{cccc}   
           e_1(T_1)  &  e_2(T_1)    \\
           e_1(t)  &  e_2(t)   \\                             
 \end{array}
 \right) \,.
$$\label{eq:det2}
That is, the 2nd order Vandermonde matrix ${\bf V}_2$ in which $T_2$ is replaced by the free variable $t$. Denote the determinant of ${\bf V}_2(T_1, t)$ by $V_2(T_1,t)$. Then the residual (\ref{eq:e2}) can be restated as
\be
r_2(t)=\frac{V_2(T_1,t)}{V_1(T_1)}\,,\label{eq:rt}
\ee
where $V_1(T_1)$ is the determinant of ${\bf V}_1$.
 
Eq.~(\ref{eq:rt}) shows that the problem of maximizing $|r_2(t)|$ (Step 6 of the algorithm) is equivalent to that one of maximizing $|V_2(T_1,t)|$, since $|V_1(T_1)|$ is a positive constant. Thus, maximizing the absolute value of the residual is equivalent to making ${\bf V}_2$ ``as invertible as possible''. 
The observation just made for the first iteration of the EIM is general, as made precise in the following theorem (see the Appendix for the proof).  \\

\begin{theorem}\label{theorem}
Define $\textup{det}({\bf V}_j):=V_j(T_1,\ldots,T_j)$. Then, the residual $r_j(t)$ computed in the $(j-1)$-iteration of the EIM-loop  satisfies 

\be\label{eq:prop}
r_j(t)=\frac{V_j(T_1,\ldots, T_{j-1},t)}{V_{j-1}(T_1,\ldots,T_{j-1})}\quad j=2,3, \ldots n \, . 
\ee
\blank
In consequence, once $T_j$ is chosen in Step 6, the residual at this node becomes
\be
r_j(T_j)=\frac{V_j(T_1,\ldots, T_{j-1},T_j)}{V_{j-1}(T_1,\ldots,T_{j-1})}=\frac{\text{det}({\bf V}_j)}{\text{det}({\bf V}_{j-1})}\,. \label{eq:residual}
\ee
\end{theorem}
This expression tells us that, at each $(j-1)$-step, the EIM algorithm selects a new node $T_j$ in order to maximize the absolute value of the determinant of ${\bf V}_j$. It is in this sense that we refer to as making the V-matrix as invertible as possible. 

To our knowledge this is a new result in the literature, and makes precise what exactly the EIM optimizes for at each iteration: it is not its conditioning or accuracy, as it is sometimes thought, but the invertibility of the associated Vandermonde matrix.
%------------------------------------------------------------------------------------------------
\subsection{Conditioning} \label{sec:cond} 
%------------------------------------------------------------------------------------------
Besides the requirement of being non-singular, for all practical purposes the V-matrix needs to be well-conditioned to any order $n$. A matrix is well-conditioned when it does not amplify small errors into large ones; in this case, the inversion of the V-matrix itself. 

In order to define a suitable quantity to measure matrix conditioning, suppose that ${\bf V}_n$  is perturbed as 
$$
{\bf V}_n  \rightarrow {\bf V}_n  + \delta {\bf V}_n\,  .  
$$
This perturbation can be due, for example, to numerical errors, which are always present in practice.
Given any matrix norm satisfying the submultiplicative property $||AB|| \leq ||A|| ||B||$ --in particular, any induced matrix norm--, one can show that
\be
\frac{ \| ({\bf V}_n + \delta {\bf V}_n)^{-1} - {\bf V}_n^{-1} \| }{ \| {\bf V}_n^{-1} \|} \leq \ka ({\bf V}_n) \frac{\| \delta{\bf V}_n \|}{ \| {\bf V}_n \|}+\mathcal{O}(\|\delta{\bf V}_n\|^2) \,, \label{eq:cond}
\ee
where the condition number of  ${\bf V}_n$ has been introduced as  
\be
\ka_n:= \ka({\bf V}_n) :=||{\bf V}_n||\,||{\bf V}_n^{-1}||\, .\label{eq:kappa1}
\ee
Note that it always satisfies $\ka({\bf V}_n) = ||{\bf V}_n||\,||{\bf V}_n^{-1}||\geq ||{\bf V}_n{\bf V}_n^{-1}|| = 1$.

The computation of the condition number not only depends on the problem but also on the norm used. Throughout this article we always use the 2--norm
$$
\|  \cdot \| := \|  \cdot \|_2\,.
$$
For briefness, we omit from hereon the $2$- subscript. 

The inequality (\ref{eq:cond}) is a standard and well known result in linear algebra. Next we turn explicitly to the conditioning of the EIM. For this we use the explicit form of the inverse of a matrix in terms of its adjugate --the transpose of the cofactor matrix--
$$
{\bf V}_n^{-1}=\frac{\text{adj}({\bf V}_n)}{|\text{det}({\bf V}_n)|}\,,
$$
to rewrite $\ka_n$ as
\be\label{eq:cond2}
\ka_n =\frac{||{\bf V}_n||\,||\text{adj}({\bf V}_n)||}{|\text{det}({\bf V}_n)|}\,.
\ee

The V-matrix is well-conditioned when its condition number $\ka({\bf V}_n)$ is ``small''. Empirically, Vandermonde-like matrices are generally ill-conditioned \cite{Gautschi2012HowU,BOYD2011443} except for some particular cases \cite{Pan2015HowBA}. However, this does not preclude  achieving stability for some matrix inversion algorithms, for example, Bj\"orck-Pereyra-type algorithms~\cite{bjorck1970,Higham1987ErrorAO,bjorck2008}. 

We return to Equation~(\ref{eq:rt}). Maximizing the absolute value of the residual $r_2(t)$ also maximizes the {\em denominator} of the condition number
\be
\ka ({\bf V}_2(T_1,t))=\frac{||{\bf V}_2(T_1,t)||\,||\text{adj}({\bf V}_2(T_1,t))||}{|V_2(T_1,t)|}\,. \label{eq:ka3}
\ee
In other words, the EIM does {\em not} find the new node $T_2$ which minimizes $\ka ({\bf V}_2(T_1,t))$, since the numerator of  (\ref{eq:ka3}) is not controlled. This observation holds for all iterations, so we restate (\ref{eq:ka3}) in its general form:
\be
\ka ({\bf V}_j(T_1,\ldots,T_{j-1},t))=\frac{||{\bf V}_j(T_1,\ldots,T_{j-1},t)||\,||\text{adj}({\bf V}_j(T_1,\ldots,T_{j-1},t))||}{|V_j(T_1,\ldots,T_{j-1},t)|}\,. \label{eq:full-ka}
\ee

So far we have shown (Theorem \ref{theorem}) that the EIM optimizes the invertibility of the V-matrix and, as we have just pointed out explicitly, not its conditioning (or accuracy, as discussed below) --as is sometimes thought. However, it does a partial job at minimizing $\ka$ (in the sense of maximizing the denominator of (\ref{eq:full-ka}))  with low computational cost. In Section \ref{sec:numerics} we discuss how this partial optimization of $\ka$ compares to the full one; that is, the minimization of the {\em full} expression (\ref{eq:full-ka}), as opposed to the maximization of its denominator only.
%------------------------------------------------------------------------------------------------
\subsection{Accuracy} \label{sec:acc}
%------------------------------------------------------------------------------------------------
We have discussed existence (Section \ref{sec:exist}) and stability/well conditioning (Section \ref{sec:cond}) of the EIM. Here we turn into its accuracy, again in the 2--norm. 

As before, assume the existence of a parameterized family of complex waveforms $h_\lam(t)$  for $t\in[0, T]$. Consider a basis $\{ e_i \}_{i=1}^n$, the generated subspace $W_n:= {\tt Span}\{ e_i \}_{i=1}^n$ and a (possibly weighted) scalar product $\langle \, . \,, .  \rangle$. The optimal linear approximation of the form
$$
 \sum_{i=1}^n c_i e_i(t)\,
$$
to $h_\lam$ is, in the least squares (LS) sense, the orthogonal projection of $h_\lam$ onto $W_n$. If we assume the basis is orthonormal, the projection becomes the familiar expression
\be\label{eq:proj}
{\cal P}_n h_\lam(t) = \sum_{i=1}^n \langle e_i,h_\lam \rangle e_i(t).  
\ee
Since (\ref{eq:proj}) is the optimal solution to the LS problem, an interpolant can only have equal or larger errors. To quantify the difference in accuracy between both schemes, the Lebesgue constant $\Lambda_n$ is usually introduced\footnote{The Lebesgue constant is usually defined in the infinity norm. In this article we define it in the 2--norm.}. First, define the discrete norm of a waveform $h_{\lam}$ as
$$
||h_{\lam}||_d^2:=\sum_{i=1}^L \bar{h_{\lam}}(t_i)h_{\lam}(t_i)\Delta t\,,
$$
where $\bar{h}$ denotes complex conjugation. The following inequality holds \cite{Antil2013TwoStepGA}
\be
\| h_\lam - {\cal I}_n [h_\lam ] \|_d \leq \Lambda_n \| h_\lam - {\cal P}_n h_\lam\|_d
 \quad\text{where}\quad
\Lambda_n := \| {\cal I}_n \|_d\,. \label{eq:lambda}
\ee
When the basis is unitary (${\bf U}{\bf U}^\dagger={\bf 1}$, with the columns of ${\bf U}$ storing the basis vectors $e_i$), the following equality is valid in the 2--norm \cite{Antil2013TwoStepGA}
\be
\Lambda_n = \|   {\bf V}^{-1}_n\|  \,. \label{eq:V}
\ee

If the approximation by projection has fast convergence, the behavior of $\Lambda_n$ with $n$ may decide how  good is the approximation by interpolation. So, it is desirable to have a slow/controlled growth of $\Lambda_n$ in order to hopefully mimic the projection approximation and then achieve fast convergence rates for the interpolation approximation as well. This is desirable, in particular, in the case of parameterized problems with regularity with respect to parameter variation as it happens for GWs: one expects (and observes), in the context of the RBM, spectral-type convergence of the Kolmogorov n-width and, as consequence, of the error decay when using greedy or POD approaches to find the basis~\cite{cohen2010,binev2011,devore2012,maday2012,rozza2013,haasdonk2013}. 
 
Notice that the Lebesgue constant stated as in (\ref{eq:V}) is directly related to the condition number $\ka({\bf V}_n)$ through \
\be
\ka({\bf V}_n)=\|{\bf V}_n\| \Lambda_n\,. \label{eq:lam}
\ee
So, controlling the conditioning of the problem and its accuracy are (partially) related. But one sees from (\ref{eq:lam}) that by optimizing $\ka$ one does not necessarily optimize for $\Lambda$, and viceversa. 

The observations of this section are generic (not specific to the EIM) and well known. Our point here is two-fold: i) to emphasize that in the same way the EIM does not optimize its conditioning it does not optimize its Lebesgue constant either, at least not in any direct way. This motivates the next Section, where we introduce two new variations of the EIM which do not optimize the invertibility of the associated Vandermonde matrix but, instead, explicitly maximize its conditioning and minimize its Lebesgue constant, respectively. 
 
%------------------------------------------------------------------------------------------------
\subsection{On Optimization of Conditioning and Accuracy} \label{sec:opt} 
%------------------------------------------------------------------------------------------------

One might wonder how the conditioning and accuracy of the EIM might improve if one replaces the optimization criteria of the EIM algorithm (step 6) by minimizing the  condition number of  ${\bf V}_j(T_1,\ldots,T_{j-1},t)$ or the Lebesgue constant at each step instead of maximizing the determinant of the V-matrix. For this purpose, we introduce two new variations of the EIM. Instead of choosing 
$$
T_j = \text{argmax}_t |r_j| \,, 
$$
with $r_j$ as in Eq.~(\ref{eq:prop}) we introduce the following two variations:
\begin{enumerate}
\item EIM-$\ka$: 
\be
T_j = \text{argmin}_t \, \ka ( {\bf V}_j(T_1,\ldots,T_{j-1},t) )  \label{eq:eim-ka}
\ee
\item EIM-$\La$: 
\be
T_j = \text{argmin}_t \, \Lambda_j(t)\quad\text{with}\quad\Lambda_j(t):=\|{\bf V}^{-1}_j(T_1,\ldots,T_{j-1},t)\| \, .  \label{eq:eim-lam}
\ee
\end{enumerate} 
That is, at each iteration, EIM-$\ka$ and  EIM-$\La$ fully minimize the condition number of the V-matrix and Lebesgue constant of the interpolant, respectively. 

A natural question that arises is how $\ka$ and $\Lambda_n$ behave with the number of basis $n$, and how they compare with the original EIM. An analytical treatment of this is beyond the scope of this work so we only present numerical results corresponding to GWs emitted by the collision of two black holes. Similar approaches for the optimization of the condition number and the the determinant of the Vandermonde matrix can be found in \cite{Guo2018WeightedAF, Shin2016NonadaptiveQP}, though they rely on polynomial approximation (not reduced basis as in this work).

%------------------------------------------------------------------------------------------------
\section{Numerical experiments} \label{sec:numerics}
%------------------------------------------------------------------------------------------------
We present numerical experiments using the EIM and our variations presented above. The study case is that one of the gravitational waves emitted by the collision of two nonspinning black holes in initial quasi-circular orbit for mass ratios $q = m_1/m_2\in [1,10]$, and time $t\in[-2750,100]$M, where M is the total mass of the system. As is usual in the field, the waveforms have been aligned so that $t=0$ corresponds to the peak of the amplitudes (roughly speaking, the time of coalescence). This corresponds to the surrogate model ${\tt SpEC\_q1\_10\_NoSpin}$~\cite{PhysRevLett.115.121102} available in~\cite{surr-catalog} and evaluated with the GWSurrogate Python package~\cite{gwsurr}.

For building the reduced basis we used the Arby Python package~\cite{arby}. Arby is a new, open source, fully data-driven module for building reduced bases, empirical interpolants and surrogate models from training data, written from scratch by one of us (AV), which will be described elsewhere. We implemented the EIM-variations EIM-$\ka$ and EIM-$\La$ in Python language. For benchmarking, we used the integration tools of Arby to perform the computation of interpolation errors of the different algorithms.

We define the maximum interpolation and projection errors in parameter space at fixed $n$ as $\tilde\sigma_n$ and $\sigma_n$, respectively. With this notation, Eq.~(\ref{eq:lambda}) reads
$$
\tilde\sigma_n\leq\La_n\sigma_n\,.
$$
The training space was filled with $1001$ equispaced complex waveforms 

$$
h(q,t):=h_+(q,t)+i\,h_{\times}(q,t)
$$
for $(q,t)\in[1,10]\times[-2750,100]$M, and time step $\Delta t=0.1$M. We recall that $q:= m_1/m_2$ is the mass ratio of the black holes and $h_+$ and $h_{\times}$ are the two transversal polarizations of a gravitational wave. This training set resulted in $19$ basis elements with squared greedy error $\sigma^2=1.202\times10^{-11}$ \footnote{We set the tolerance error to $\sim 10^{-11} - 10^{-12}$ before reaching double precision so as to avoid roundoff issues.}. 

\begin{figure}[H]
\begin{center}
\includegraphics[width=0.45\columnwidth]{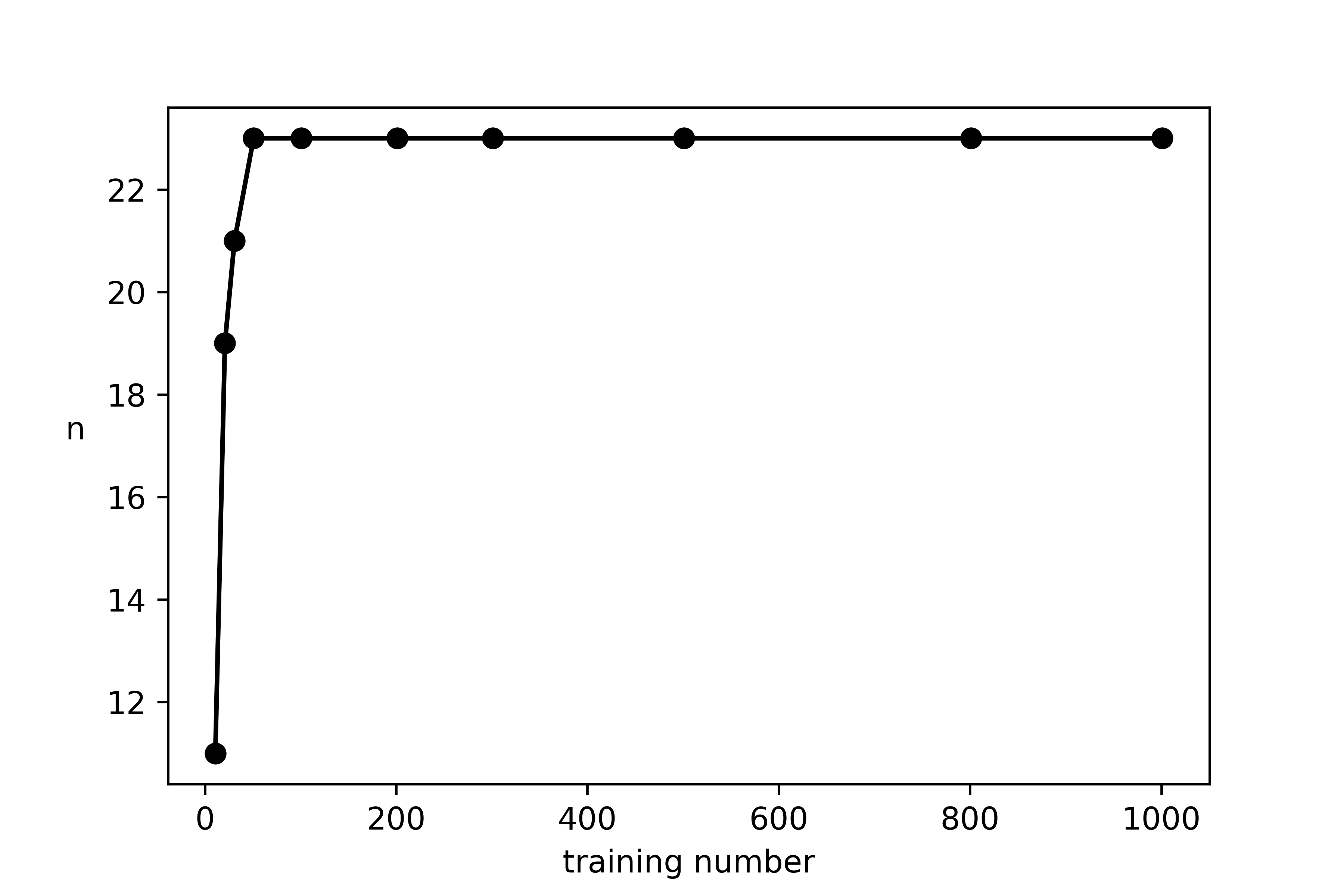}
\caption{Number of basis elements $n$ vs. the number of training elements associated to equispaced training spaces up to $1001$. The step function behavior of the curve is due to the fact that the underlying model to reduce is already spanned by a finite-dimensional space of waveforms.}
\label{fig:number}
\end{center}
\end{figure}

Note that the underlying model used to populate the training space is already a surrogate one, built with the RBM and the EIM~\cite{PhysRevLett.115.121102}. This means that the entire model belongs to a finite dimensional function space spanned by some basis. This is illustrated in Figure \ref{fig:number}. For a greedy tolerance error of $10^{-16}$, the number of basis elements $n$ needed to construct a reduced basis for a given training set rapidly reaches the value $n=23$. The step function behavior of Figure  \ref{fig:number} can be explained as follows: if we set the $23$ greedy waveforms as our initial training set and next we populate it with arbitrary waveforms from the surrogate model, subsequent training spaces will represent linearly dependent waveforms of the $23$ initial ones and the curve will become an exact step function, as it nearly happens to be the case. 
The three algorithms, EIM, EIM-$\La$ and EIM-$\ka$ choose qualitatively similar nodes; in all cases with clustering near the boundaries of the time interval, as expected \cite{hesthaven2007spectral}. Given that the selection of the nodes is very similar in the three cases, one would not expect (and we do not find) major differences between them. The underlying, rigorous, reason for this is unknown, since their optimization criteria are very different. A deeper mathematical understanding is needed but is beyond the scope of this work. 

\begin{figure}[H]
\begin{center}
\includegraphics[width=0.45\columnwidth]{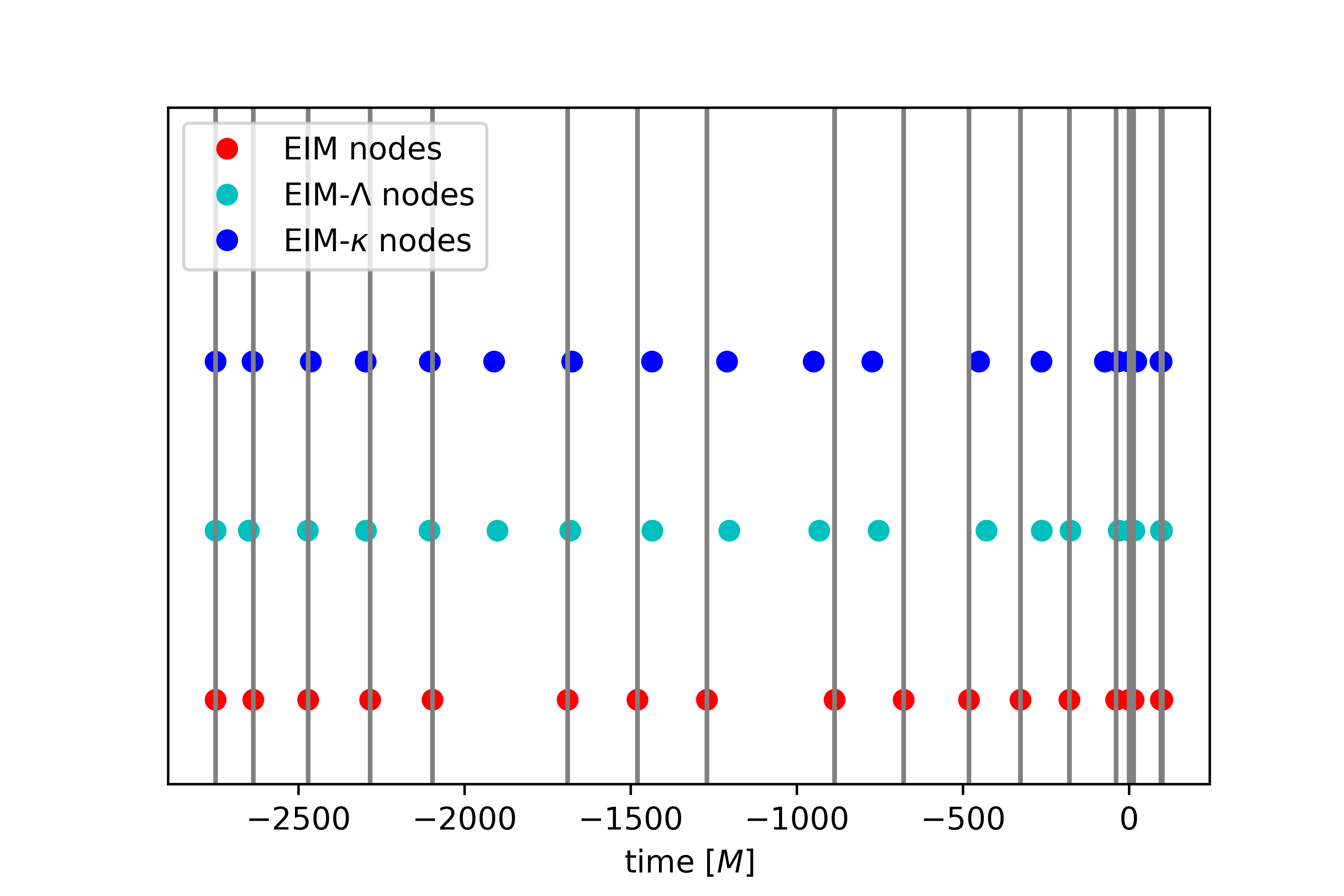}
\includegraphics[width=0.45\columnwidth]{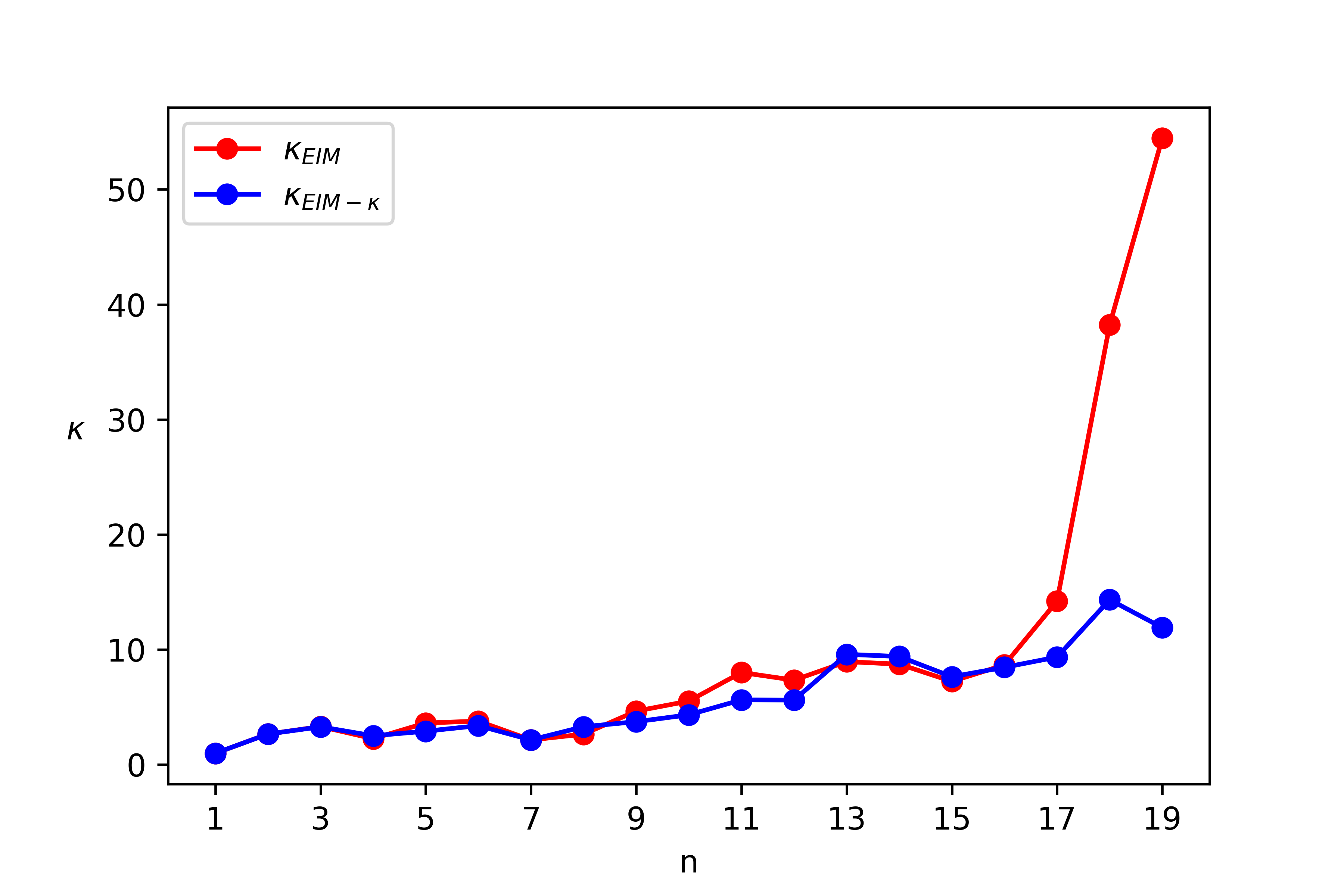}
\caption{{\bf Left}: Distribution of the $19$ interpolation nodes for the three algorithms with vertical lines centered at each EIM node. {\bf Right}: Condition number dependence with basis number $n$ for both EIM and EIM-$\ka$.  EIM-$\ka$ does show better behavior of $\ka$, as expected by construction, but only significantly at very high resolutions.}
\label{fig:gw1}
\end{center}
\end{figure}

Figure~\ref{fig:gw1} shows these distributions alongside the curves comparing the behavior of the condition number $\ka$ with the number of basis elements $n$ for both the EIM and EIM-$\ka$. One can notice that, as expected, the conditioning of the latter does improve, but only significantly at high resolutions.  

\begin{figure}[H]
\begin{center}
\includegraphics[width=0.45\columnwidth]{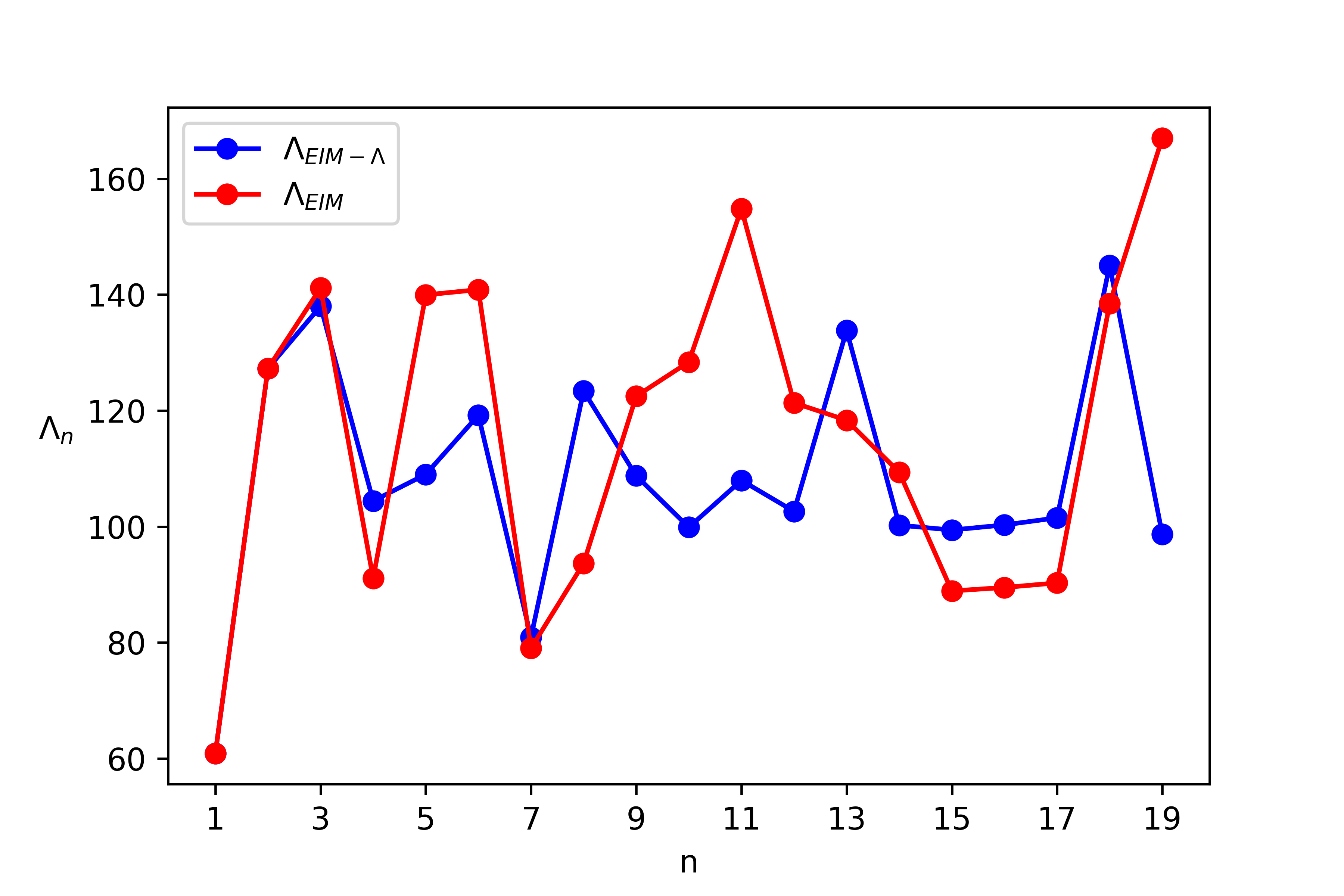}
\includegraphics[width=0.45\columnwidth]{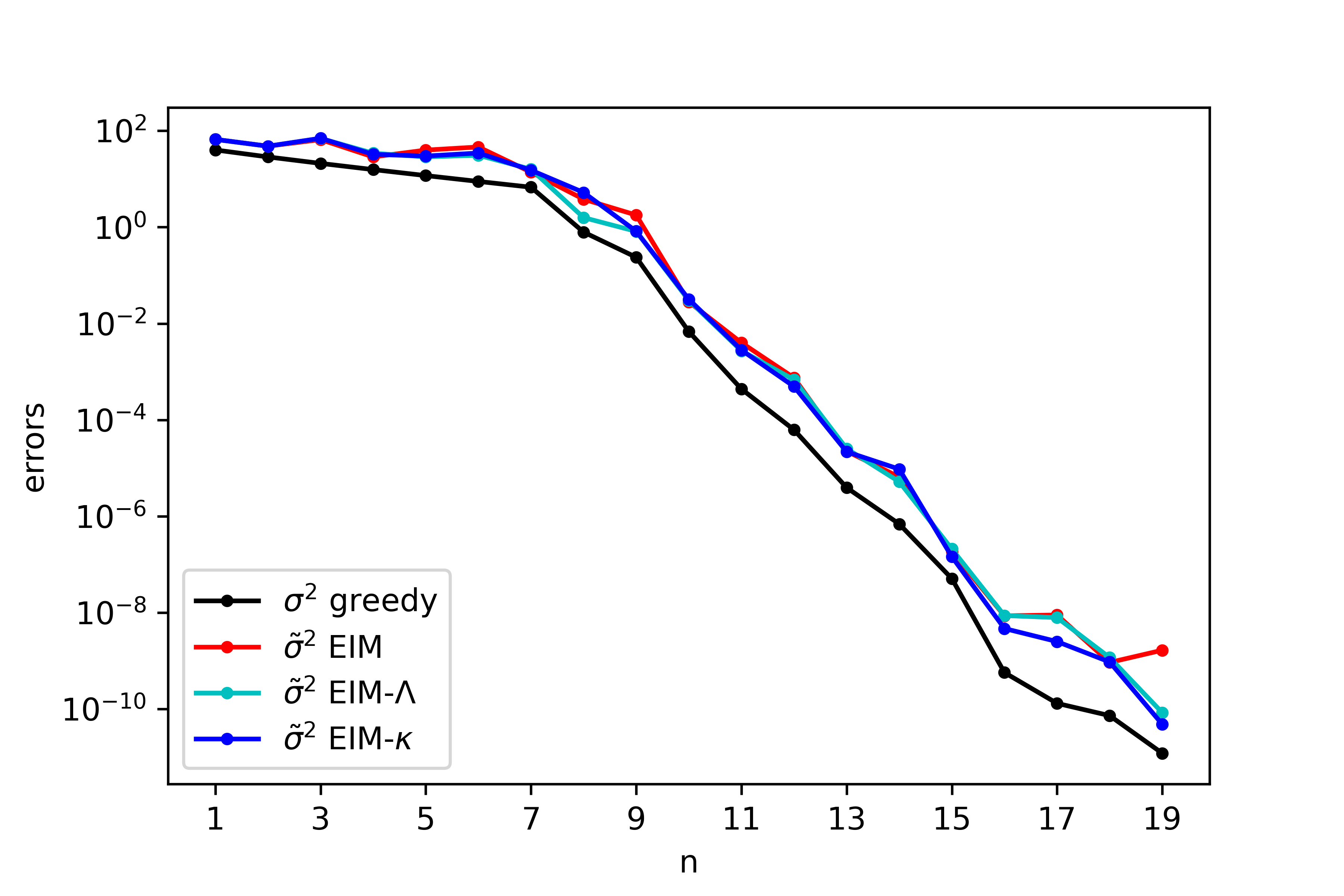}
\caption{{\bf Left}: Lebesgue constant $\La_n$ dependence with $n$ for both EIM and EIM-$\La$. {\bf Right}: In-sample validation for the three algorithms. Dots on the right panel correspond to the maximum squared interpolation errors $\tilde\sigma_n^2$ and squared greedy errors $\sigma_n^2$. Note that the greedy errors bounds from below the interpolating errors as expected, since the LS problem is optimal in the 2--norm.}
\label{fig:gw2}
\end{center}
\end{figure}

Figure \ref{fig:gw2} shows the dependence of $\La_n$ with $n$ for both EIM and EIM-$\La$. Here there is no apparent improvement; in fact, the EIM performs sometimes slightly better, sometimes slightly worse, than EIM-$\La$. This apparent contradiction is subtle and requires some thought. It is ultimately a consequence of the nested character of the algorithms: at each step, the next node relies on the previous ones, so, there is no a step-to-step updating of nodes relative to an actual iteration, as would happen in a global optimization approach. More explicitly, given nodes 
$$
\{T_1, \ldots T_{j-1} \}\,
$$
EIM-$\La$ chooses the next node as in  (\ref{eq:eim-lam}) to minimize the Lebesgue constant $\La_j$ having kept all the previous nodes fixed. This does not imply that $\La_j$ is smaller with respect to a different choice of nodes $\{T_1, \ldots T_j\}$ for all $j$, for example those of the original EIM, though we would have expected so. This is in contrast to a global optimization approach, where all nodes are chosen together; but such an approach is not only computationally expensive but also loses one of the main features of the RBM and the EIM: being nested. 

Figure \ref{fig:gw2}  also shows the maximum (over the training set) squared interpolation errors $\tilde{\sigma}_n^2$ for each $n$ corresponding to the three versions of the EIM. We also show the RBM squared greedy errors $\sigma_n^2$. The latter, as expected, is a lower bound for all the interpolation approximations. The point to notice, beyond the exponential convergence, is that all curves are very close to each other, suggesting that the EIM and our variations are nearly optimal, a rigorous proof of which is not available at the time of this writing. The reason why the errors of the EIM, as shown in Figure~\ref{fig:gw2} --and as found in practice not only in our experiments but in most if not all cases where the EIM is used--, is not monotonically decreasing is the same one we have just discussed with respect to the Lebesgue constant (besides the fact that it does not necessarily optimize for accuracy). 

\begin{figure}[H]
\begin{center}
\includegraphics[width=0.45\columnwidth]{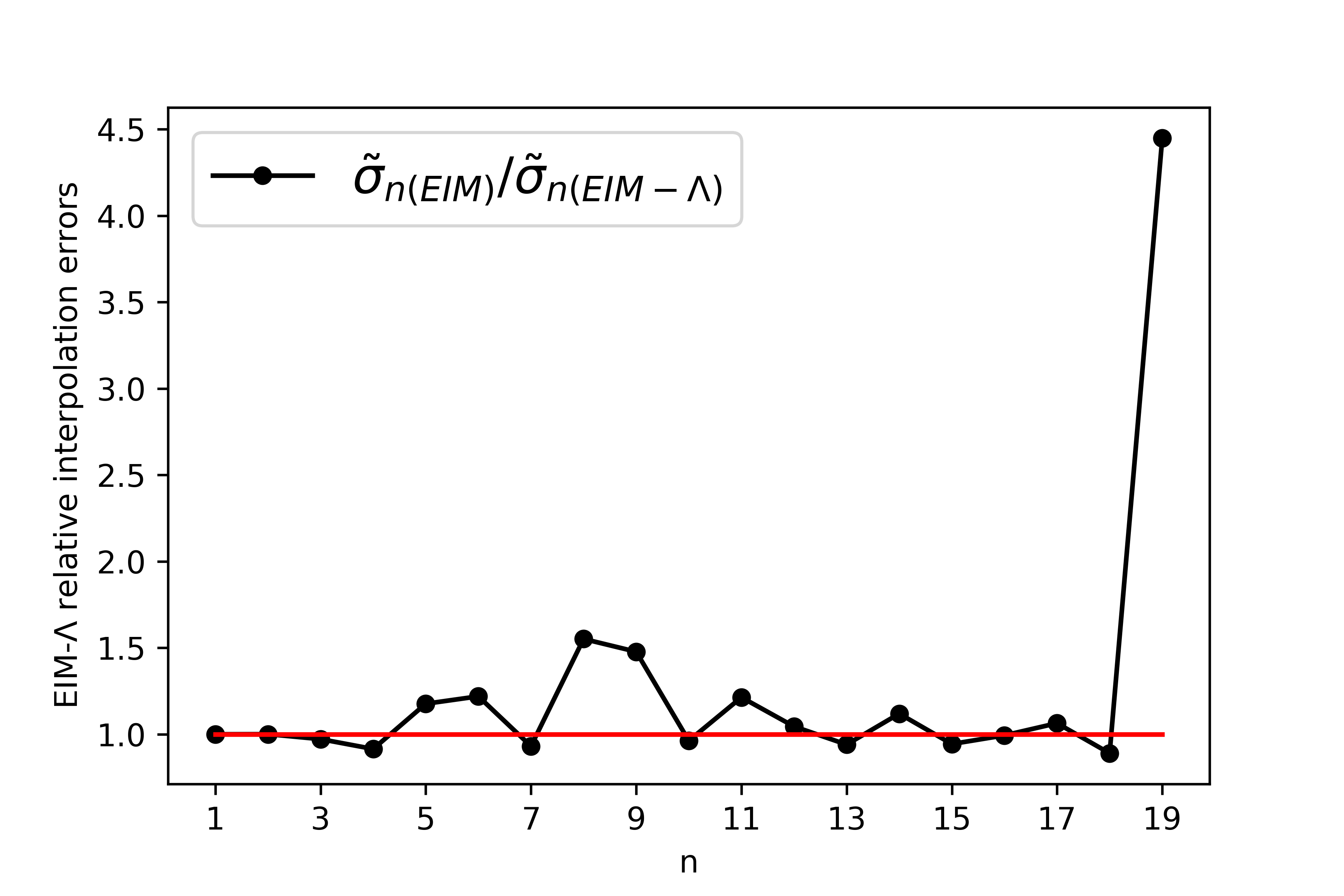}
\includegraphics[width=0.45\columnwidth]{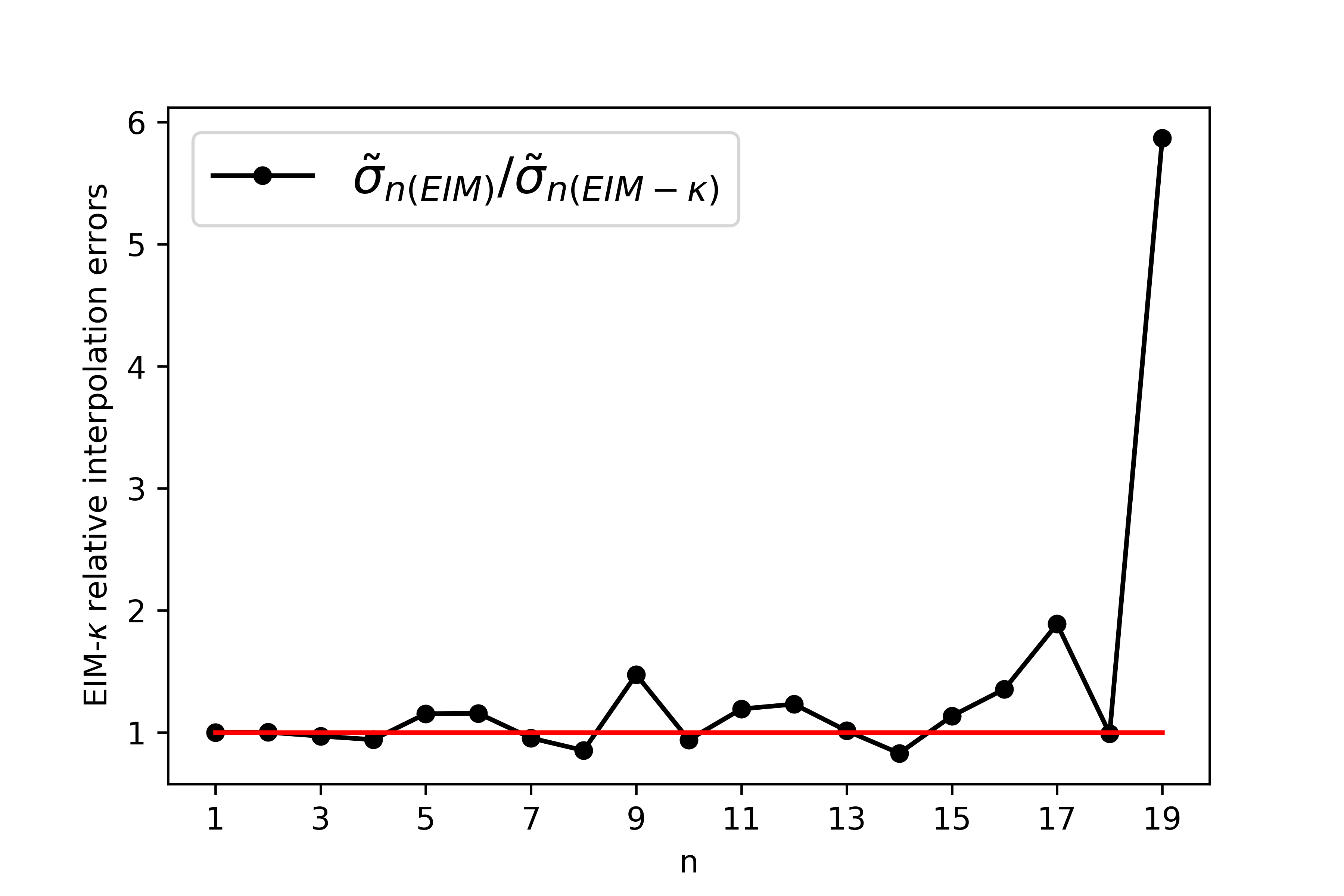}
\caption{Ratios $\tilde\sigma_{n(EIM)}/\tilde\sigma_{n(EIM-\La/\kappa)}$ between maximum interpolation errors for original EIM and the EIM$-\ka$ and EIM$-\La$ algorithms. They are qualitatively similar, except at very high resolutions.}
\label{fig:gw3}
\end{center}
\end{figure}

The Lebesgue constant provides a bound for the interpolation error in terms of the projection error. Since it is a bound, it is worthwhile looking at the actual errors themselves. Figure~\ref{fig:gw3} shows the quotients between maximum interpolation errors for EIM-$\La$ and EIM-$\ka$ when compared to the original EIM. There is an improvement on average, especially at high resolutions, but it is not significative. This can actually already be seen qualitatively  in the right panel of  Fig.~\ref{fig:gw2}. 

There is no obvious reason why on average EIM-$\ka$  performs --though marginally-- better in accuracy, since it does not optimize for it. We would have expected an improvement with EIM-$\La$, though. In hindsight the Lebesgue constant provides a bound for the error of the interpolant, and it does not need to reflect the actual error.   

%------------------------------------------------------------------------------------------------
\section{Comments} \label{sec:comments}
%------------------------------------------------------------------------------------------------

We first proved our main result: a relationship presented in Section~\ref{sec:exist},  between the EIM and the determinant of the associated Vandermonde matrix, which shows that the EIM as originally introduced and used in practice, at least in GW science, maximizes the invertibility of its associated Vandermonde matrix, not its conditioning or accuracy of the resulting interpolant, as it might be generally thought. Motivated by this result, we explored two new different optimization criteria applied to variations of the EIM, nested as well, in the context of the Reduced Basis Method.  

We optimized the EIM with respect to the condition number of the associated Vandermonde matrix and with respect to the Lebesgue constant of the resulting interpolant, both in the 2--norm, variations which we named EIM-$\ka$ and EIM-$\La$, respectively. For the original EIM applied to GWs, we have found a rapid growth of the condition number $\ka$ at high resolutions (right panel of Figure~\ref{fig:gw1}), which can be avoided through a simple modification, EIM-$\ka$. This result might be worth keeping in mind when working with high accuracy. The very slow growth of the Lebesgue constant in the 2--norm showed in these experiments is consistent with other numerical results outside GW science, see for example Figure 8 of \cite{Maday2013}. 

Our variations of the EIM have minimal changes in terms of computational cost. The online cost of the surrogate evaluation (\ref{eq:interp-sol1}, \ref{eq:interp-sol2}) is independent of the approach used to choose the interpolation nodes. The offline cost of our variations might be slightly higher, since they involve computing the condition number and Lebesgue constant at each iteration, see Eqs.~(\ref{eq:eim-ka} and (\ref{eq:eim-lam}). However, since at least in our domain of interest the bases are quite sparse, even in higher dimensions than those considered here, this offline computational cost is negligible. Furthermore, being offline, the construction of the interpolant has to be done only once, similarly to the construction of the reduced basis, and the orders of magnitude speedup of the resulting surrogates compared to supercomputer simulations of the Einstein equations justify those offline costs.  

Even though these improvements might be subtle in practice, we have shown what exactly the EIM optimizes for, which appears to be contrary to common perception. As of this writing and up to our knowledge there are no rigorous results or theory for the exact reasons, or under which conditions, the EIM as originally introduced shows well conditioning and, in practice, a controlled growth of the Lebesgue constant, since it does not optimize for either of them. 

A deeper theoretical understanding in the context of the EIM appears to be needed, since it seems that there is a subtle relationship between existence of the interpolant, conditioning of the related Vandermonde matrix and the Lebesgue constant of EIM interpolants. Namely, maximizing (\ref{eq:residual}), minimizing (\ref{eq:full-ka}) or (\ref{eq:V}) at {\em each} iteration, and how the different full nested approaches compare with each other for any $n$. Many of these points are in fact active areas of research in approximation theory, see for example \cite{maday2015} for the Generalized Empirical Interpolation Method in the context of Banach spaces. Since Reduced Order Modeling has proved to be a pivotal and relatively new field in GW science, enabling real time evaluations of GWs from compact binaries, to quasi real time parameter estimation (for astrophysical purposes, in the form of rapid followup of electromagnetic counterparts) through Focused Reduced Order Quadratures \cite{PhysRevD.102.104020}, our goal in this paper has been to introduce and prove a new result \ref{theorem}, and make accessible to practitioners some of the technical underlying and open issues through some new analyses and numerical experiments.

----------------------------------------------------------------------------------------------
\section*{Acknowledgments}
We thank Jorge Pullin, the Horace Hearne Institute and the Center for Computation and Technology at LSU, for hospitality while part of this work was done. We also thank Scott Field, Stephen Lau, and two anonymous referees for valuable comments on a previous version of this manuscript. This work was supported partially by CONICET and EVC-CIN (Argentina).
%------------------------------------------------------------------------------------------------

\blank
\par\noindent\rule{\textwidth}{0.4pt}
\blank

\appendix

\section{Proof of Theorem \ref{theorem}} 
 Consider the nested matrix ${\bf V}_j(t)$ defined as

\begin{equation*} \label{eq:Vander-t}
  {\bf V}_j(t) := \left(  \begin{array}{cccc}   
                     				 &&& e_j(T_1)      \\
                   &{\bf V}_{j-1}    	&& e_j(T_2)      \\              
                        				&&& \vdots        \\   
  e_1(t)&e_2(t)               &\cdots & e_j(t)        \\               
             \end{array}
   \right) \,,
\end{equation*}
where ${\bf V}_{j-1}$ is the $(j-1)$-order V-matrix associated to the empirical nodes $T_1,\ldots, T_{j-1}$. That is, ${\bf V}_j(t)$ is ${\bf V}_j$ with $T_j$ replaced by $t$. Write the determinant of ${\bf V}_j(t)$ as
\be\label{eq:det1}
V_j(T_1,\ldots, T_{j-1},t)=e_j(t) \textup{det}({\bf V}_{j-1}) + e_j(T_{j-1})C_{j-1\,n}[{\bf V}_j (t)]+\ldots+e_j(T_1)C_{1\,j}[{\bf V}_j (t)]\,,
\ee
where 
$$
C_{ik}[\cdot]=(-1)^{i+k}D_{ik}[\cdot]
$$ 
denotes the $(i,k)$-cofactor of $[\cdot]$, and $D_{ik}[\cdot]$ denotes the determinant of $[\cdot]$ after eliminating its $i$-th row and $j$-th column. 

Divide (\ref{eq:det1}) by $\textup{det}({\bf V}_{j-1})=V_{j-1}(T_1,\ldots,T_{j-1})$\ft{The EIM-loop ensures, by independence of row and column vectors, that the determinant of the V-matrix is non zero to any order $j\geq 1$.}

\begin{align}\label{eq:det2}
\frac{V_j(T_1,\ldots, T_{j-1},t)}{V_{j-1}(T_1,\ldots,T_{j-1})}=e_j(t)-\frac{1}{\textup{det}({\bf V}_{j-1})}\Big\{e_j(T_{j-1})D_{j-1\,j}[{\bf V}_j (t)]+\ldots &\\
+(-1)^{j}e_j(T_1)D_{1\,j}[{\bf V}_j (t)]\Big\}&\nonumber\,.
\end{align}
Let's write the substracting term of the r.h.s. of (\ref{eq:det2}) in matrix notation,
\be\label{eq:det}
\frac{1}{\textup{det}({\bf V}_{j-1})}\begin{pmatrix}\,\,(-1)^j D_{1\,j}[{\bf V}_j (t)]&\cdots&  D_{j-1\,j}[{\bf V}_j (t)] \,\,  \end{pmatrix}\begin{pmatrix}e_j(T_1)\\ \vdots \\e_j(T_{j-1})\end{pmatrix}\,.
\ee
Now write the interpolant $\cI_{j-1}[e_j](t)$ in the same spirit,

\begin{align*}
\cI_{j-1}[e_j](t)&=\begin{pmatrix}e_1(t)&\ldots&e_{j-1}(t)\end{pmatrix}{\bf V}_{j-1}^{-1}\begin{pmatrix}e_j(T_1)\\ \vdots \\e_j(T_{j-1})\end{pmatrix}\\
&=\frac{1}{\textup{det}({\bf V}_{j-1})}\begin{pmatrix}e_1(t)&\ldots&e_{j-1}(t)\end{pmatrix}\textup{adj}({\bf V}_{j-1})\begin{pmatrix}e_j(T_1)\\ \vdots \\e_j(T_{j-1})\end{pmatrix}\,. \nonumber
\end{align*}
This is similar to expression \ref{eq:det} above. If one proves that 
\be\label{eq:1}
\begin{pmatrix}\,\,(-1)^j D_{1\,j}[{\bf V}_j (t)]&\cdots&  D_{j-1\,j}[{\bf V}_j (t)] \,\, \end{pmatrix}
\ee
is equal to  
\be\label{eq:2}
\begin{pmatrix} P_1(t)&\ldots&P_{j-1}(t) \end{pmatrix}:=\begin{pmatrix}e_1(t)&\ldots&e_{j-1}(t)\end{pmatrix}\textup{adj}({\bf V}_{j-1})\,,
\ee
the proof of formula \ref{eq:prop} will be completed. Let's look at the first component of (\ref{eq:1}):
\begin{equation*}
(-1)^j D_{1\,j}[{\bf V}_j (t)]=\sum_{i=1}^{j-1}(-1)^{i-1}e_i(t)\,D_{j-1\,i}[{\bf V}_j (t) (1|j)]\,,
\end{equation*}
where ${\bf V}_j (t)(1|j)$ stands for ${\bf V}_{j}(t)$ without its 1-th row and $j$-th column. On the other hand, the first element of (\ref{eq:2}) is equal to
\begin{align*}
P_1(t)=&\sum_{i=1}^{j-1}e_i(t)(\textup{adj}({\bf V}_{j-1}))_{i\,1}=\sum_{i=1}^{j-1}e_i(t)C_{1\,i}({\bf V}_{j-1})\\
&=\sum_{i=1}^{j-1}(-1)^{i-1} e_i(t)D_{1\,i}({\bf V}_{j-1})\,.\nonumber
\end{align*}
Finally, notice that $$D_{j-1\,i}[{\bf V}_j (t) (1|j)]=D_{1\,i}({\bf V}_{j-1})\,.$$ Therefore $$(-1)^j D_{1\,j}({\bf V}_j (t))=P_1(t)\,$$ and, in the same way, $$(-1)^{j-1-i} D_{i\,j}({\bf V}_j (t))=P_i(t)\quad \text{for}\,\, i=1,\ldots, j-1\,.$$ Therefore $\cI_{j-1}[e_j](t)$ is equal to expression \ref{eq:det}:
$$
\frac{V_j(T_1,\ldots, T_{j-1},t)}{V_{j-1}(T_1,\ldots,T_{j-1})}=e_j(t)-\cI_{j-1}[e_j](t)=r_j(t)
$$
and the Theorem is proved.

\bibliography{refs}

\end{document}